\documentclass[12 pt,amsfonts] {amsart}
\def\dom{\mathop{\mathrm{dom}}\nolimits}
\def\cod{\mathop{\mathrm{cod}}\nolimits}
\def\img{\mathop{\mathrm{img}}\nolimits}
\def\rng{\mathop{\mathrm{rng}}\nolimits}
\def\id{\mathop{\mathrm{id}}\nolimits}
\def\R{{\Bbb R}}
\title{Atlas as solution of Sincov's inequality}
\author{Petra Augustov\'a and Lubom\'{\i}r Klapka}
\begin{document}

\maketitle

\noindent{\bf Abstract:} 
We find a general solution of Sincov's inequality
$\Phi_{\alpha\beta}\circ\Phi_{\beta\gamma}\subseteq
\Phi_{\alpha\gamma}$ provided
$\Phi_{\alpha\beta}^{-1}\subseteq\Phi_{\beta\alpha},
\Phi_{\alpha\alpha}\subseteq\id_M$. Further, we prove that in the differentiable 
case we can interpret such solution as a differentiable manifold in the original 
sense of Lang. This allows to generalize the notion of atlas and transition map 
for non-differentiable and discontinuous case. 

\bigskip


\bigskip


\section*{Sincov's functional inequality}

\noindent We are going to solve functional inequalities
\begin{equation}
\label{1} \Phi_{\alpha\beta}\circ\Phi_{\beta\gamma}\subseteq
\Phi_{\alpha\gamma},\quad
\Phi_{\alpha\beta}^{-1}\subseteq\Phi_{\beta\alpha},\quad
\Phi_{\alpha\alpha}\subseteq\id_M
\end{equation}
for an unknown map 
\begin{equation}
\label{2} \Phi\colon I\times
I\ni(\alpha,\beta)\mapsto\Phi_{\alpha\beta}\in2^{M\times M}
\end{equation}
where $I$ and $M$ are sets, $\circ$ denotes composition of binary relations and 
the exponent $-1$ stands for inverse binary relation.

\bigskip

\subsection*{Motivation example: Dependence of the solution of an ordinary differential equation on initial condition.} 

Consider an ordinary differential equation (ODE)
$\dot x=f(\tau,x)$ where $x$ is a dependent variable, $\tau$ an independent variable, 
$f$ a differentiable map of class  $C^1$ with domain $\dom f\subseteq\R^2$ being an open set and
codomain $\cod f=\R$. 
Since any Cauchy problem $x(\alpha)=a$
where $(\alpha,a)\in\dom f$, has in this case a single maximal solution, 
it gives a map $x=F(\tau,\alpha,a)$, where $\dom
F\subseteq\R^3$ is an open set and $\cod F=\R$ (see Theorem 14
in section 23  \cite{1}). This map belongs to the differentiability class 
$C^1$ (see Corollary 4 in section 7 \cite{2}). If variables $\tau$,
$\alpha$, $\beta$ and $a$ satisfy conditions $(\tau,\alpha,a)\in\dom
F$ and $(\beta,\alpha,a)\in\dom F$, then for the map $F$ the following incidences
hold:

1) Since it is a solution of Cauchy problem, we have 
$$F(\alpha,\alpha,a)=a.$$

2) Since Cauchy problem has a unique solution, we also have 
$$F(\tau,\alpha,a)=F(\tau,\beta,F(\beta,\alpha,a)).$$

Those incidences cannot be considered in general as functional equations of an unknown map $F$ because they incorporate implications on the domain of $F$,
especially we have for 1) $(\tau,\alpha,a)\in\dom
F\Rightarrow(\alpha,\alpha,a)\in\dom F$ and we have for 2) $(\tau,\alpha,a)\in\dom F\,\,\wedge\,\,(\beta,\alpha,a)\in\dom
F\,\,\Rightarrow\,\,(\tau,\beta,F(\beta,\alpha,a))\in\dom F.$
Thus, unless we know $F$ we don't know all possible values of variables $\tau$, $\alpha$, $\beta$ and $a$. 
Hence, the problem is about searching for $F$ and $\dom F$ at the same time. 

This problem can be rigorously treated by means of binary relations provided that by 
{\it binary relation\/} we mean the set of ordered pairs (see, e.g., Definition 1. v \S3.1
\cite{3} or 6.8 Definition in Chapter I \cite{4}). In the sequel, we consider binary relations precisely in this sense.  
By {\it composition\/} $\rho\circ\sigma$ of binary relations $\rho$ and $\sigma$ we mean a binary relation 
$\rho\circ \sigma=\{(a,b)|\,\,\exists c\colon (c,b)\in \rho\,\,\wedge\,\,(a,c)\in \sigma\}$, by {\it inverse\/}
$\rho^{-1}$ of a binary relation $\rho$ a binary relation 
$\rho^{-1}=\{(a,b)| (b,a)\in \rho\}$, by {\it domain\/} $\dom
\rho$ of a binary relation $\rho$ the set $\dom\rho=\{a|\,\,\exists
b\colon(a,b)\in \rho\}$ and by {\it range\/} $\rng \rho$ of a binary
relation $\rho$ the set $\rng \rho=\{b|\,\,\exists a\colon(a,b)\in
\rho\}$.

We introduce a system
$\{\Phi_{\alpha\beta}\}_{(\alpha,\beta)\in\R^2}$ of binary relations
$\Phi_{\alpha\beta}$ defined by  
$$
(b,a)\in\Phi_{\alpha\beta}\quad\Leftrightarrow\quad
a=F(\alpha,\beta,b).
$$
Then the inequalities $(\ref{1})$ are for any choice of dependent variables $\alpha,\beta, \gamma\in\R$ 
equivalent to the incidences for $F$ and implications for $\dom F$ that are considered in this example. 
Hence, we can treat them as functional inequalities for the map $\Phi$ defined by $(\ref{2})$ where $I=M=\R$.

Once we find the map $\Phi$, we can return to the map $F$. 
Then the domain of $F$ is given by
$$
b\in\dom\Phi_{\alpha\beta}\quad\Leftrightarrow\quad
(\alpha,\beta,b)\in\dom F.
$$
Moreover, we can return also to the original differential equation. 
To do that, we differentiate the incidence 2) with respect to $\tau$ and we set 
$\beta=\tau$, then $\dot x=f(\tau,x)$, where
$$
f(\tau,x)=\left.\frac{\partial
F(\tau,\beta,x)}{\partial\tau}\right|_{\beta=\tau}.
$$

\section*{General solution}

If we have equalities in $(\ref{1})$, then we deal with a {\it modified 
Sincov's equation\/} $\Phi_{\alpha\beta}\circ
\Phi_{\beta\gamma}=\Phi_{\alpha\gamma}$ with known solution 
(see, e.g., Theorem 2 v 8.1.3 \cite{5}). In our work we consider another 
modification of the modified Sincov's equation that we call 
{\it  modified Sincov's inequality\/} $\Phi_{\alpha\beta}\circ
\Phi_{\beta\gamma}\subseteq\Phi_{\alpha\gamma}$ and we find its solution 
under condition $\Phi_{\alpha\beta}^{-1}\subseteq
\Phi_{\beta\alpha}$ and $\Phi_{\alpha\alpha}\subseteq\id_M$.

Before we formulate the result, we need to define an {\it injective
binary relation\/} or shortly {\it injection\/}. By that we understand 
a  binary relation $\rho$ with the property $(b_1,a)\in
\rho\,\wedge(b_2,a)\in \rho\,\Rightarrow b_1=b_2$. This is  
connected with the notion of co-injection. A {\it co-injective binary relation\/} 
or shortly {\it co-injection\/} is such a binary relation $\rho$ that
$\rho^{-1}$ is an injection. We can easily see that the injectivity criterion 
of  the relation $\rho$ is $\rho^{-1}\circ\rho=\id_{\dom\rho}$ and its 
co-injection criterion $\rho\circ\rho^{-1}=\id_{\rng\rho}$.

\bigskip

\noindent{\bf Theorem 1:} {\it The map $(\ref{2})$ is a solution of 
functional inequalities $(\ref{1})$ if and only if there exist such 
a system $\{\varphi_\alpha\}_{\alpha\in I}$ of injective co-injections
$\varphi_\alpha$ such that
\begin{equation}
\label{3} \Phi\colon I\times I\ni(\alpha,\beta)
\mapsto\varphi_\alpha\circ\varphi^{-1}_\beta\in 2^{M\times M}.
\end{equation}
\/}

\noindent{\bf Proof:} We assume that the map $(\ref{2})$ is 
a solution of functional inequalities $(\ref{1})$ and we define a binary 
relation $\psi$ by
$$
((\beta,b),(\alpha,a))\in\psi\,\,\Leftrightarrow\,\,(b,a)\in
\Phi_{\alpha\beta}.
$$
Since $\Phi_{\alpha\beta}\circ\Phi_{\beta\gamma}\subseteq
\Phi_{\alpha\gamma}$, the relation is transitive. Since further
$\Phi_{\alpha\beta}^{-1}\subseteq\Phi_{\beta\alpha}$, this relation
is also symmetric. Any binary relation that is transitive and symmetric, 
is an equivalence on its range, hence we can decompose 
$\psi$  by means of the corresponding quotient projection 
$\pi\colon\rng\psi\,\to\,\rng\psi/ \psi$ into a composition 
$\psi=\pi^{-1}\circ\pi$. Now, we define for any $\alpha\in I$
a binary relation $\varphi_\alpha$ by 
$$
(z,a)\in\varphi_\alpha\,\,\Leftrightarrow\,\,((\alpha,a),z)\in\pi .
$$
Hence, we have $(b,a)\in\varphi_\alpha
\circ\varphi_\beta^{-1}\,\Leftrightarrow\,((\beta,b),(\alpha,a))\in\psi$,
and therefore (\ref{3}) is proved. Since the projection $\pi$ is a 
co-injection, then any binary relation $\varphi_\alpha$  is injective.
And since $\varphi_\alpha\circ\varphi^{-1}_\alpha=
\Phi_{\alpha\alpha}\subseteq\id_M$, any relation 
$\varphi_\alpha$ is a co-injection. This proves one implication.

The inverse implication follows from substituting  $(\ref{2})$, $(\ref{3})$ into 
$(\ref{1})$ and using equalities
$$
\varphi^{-1}_\alpha\circ\varphi_\alpha=\id_{\dom\varphi_\alpha},\quad
\varphi_\alpha\circ\varphi^{-1}_\alpha=\id_{\rng\varphi_\alpha},\quad
\rng\varphi_\alpha\subseteq M,
$$
which express that $\{\varphi_\alpha\}_{\alpha\in I}$ is a system
of injective co-injections satisfying $(\ref{3})$. $\Box$

\section*{Properties of injective co-injections}

Injective co-injections play an important role in finding solutions of 
functional inequalities $(\ref{1})$. Therefore,  we will first recapitulate
some of the most important corresponding properties 
that are easy to verify from definitions.

First of all, injective co-injections form a subalgebra of the algebra of 
binary relations together with composition and inverse operation.
In the context of Theorem 1, this means that 
not only the relation $\varphi_\alpha$ but also all relations
$\Phi_{\alpha\beta}$ solving the inequalities $(\ref{1})$ are
injective co-injections.

To link it with known results it is also important to realize 
connections with maps. In particular,  co-injections are closely 
related to surjective maps and injective co-injections to bijective maps.
We will describe it in detail in the following. 

For any binary relation $\rho$ holds
$\rho\subseteq\dom\rho\times\rng\rho$. By definition of the domain, 
for any $b\in\dom\rho$ there exists $a$ such that $(b,a)\in\rho$. If the 
relation $\rho$ is co-injective, there is exactly one such $a$.

This allows to define a map  $\tilde{\rho}\colon\dom\rho\to\rng\rho$ 
that assigns to any $b\in\dom\rho$ a unique $a\in\rng\rho$ such that 
\begin{equation}
\label{6} (b,a)\in\rho\,\Leftrightarrow\,\tilde{\rho}(b)=a.
\end{equation}
By definition of the range, this map is surjective. On the other hand, 
any surjective map 
$\tilde{\rho}\colon\dom\tilde{\rho}\to\cod\tilde{\rho}$ 
determines a co-injective binary relation 
$\rho\subseteq\dom \tilde{\rho}\times\cod\tilde{\rho}$ satisfying
$(\ref{6})$.

Since the surjection $\tilde{\rho}$ is fully described by the 
co-injection $\rho$, without loss of accuracy in the following,  
we can identify the co-injection $\rho$ with the surjection
$\dom\rho$ on $\rng\rho$ and use for co-injections the notation 
$\rho(b)=a$ as alternative to $(b,a)\in\rho$.

Moreover, if the co-injection $\rho$ is injective, then it is 
bijective as a map in the previously described sense. 

The next theorem treats the relation between collections of 
injective co-injections giving by means of $(\ref{3})$ the same
solution of functional inequalities $(\ref{1})$.

\bigskip

\noindent{\bf Theorem 2:} {\it Let $\{\varphi_\alpha\}_
{\alpha\in I}$ and $\{\bar{\varphi}_\alpha\}_{\alpha\in I}$ be
two sets of injective co-injections such that for any
$(\alpha,\beta)\in I\times I$ holds
\begin{equation}
\label{4} \varphi_\alpha\circ\varphi^{-1}_\beta=
\bar{\varphi}_\alpha\circ\bar{\varphi}^{-1}_\beta.
\end{equation}
Then there exists an injective co-injection $\omega$ such that 
for any $\alpha\in I$ 
\begin{equation}
\label{5} \varphi_\alpha=\bar{\varphi}_\alpha\circ\omega.
\end{equation}
\/}

\noindent{\bf Proof:} Let $\chi$ be a binary relation $\chi$ 
defined by 
$$
((\alpha,a),z)\in\chi\,\,\Leftrightarrow\,\,(z,a)\in\varphi_\alpha.
$$
Since the relation $\varphi_\alpha$ is injective, the relation 
$\chi$ is a co-injection. Therefore, the symmetric relation 
$\chi^{-1}\circ\chi$ is transitive and we have an equivalence on 
$\rng(\chi^{-1}\circ\chi)=\dom\chi$. The canonical decomposition 
of co-injection $\chi$ treated as a surjective map
$\chi\colon\dom\chi\to\rng\chi$ is of the form 
$\chi=\nu\circ\vartheta\circ\pi$. The surjectivity implies the identity 
of the embedding $\nu$, hence we can compute the 
projection $\pi=\vartheta^{-1}\circ\chi$, where $\vartheta$ is a 
bijection of the canonical decomposition of surjection $\chi$.

Furthermore, if we introduce a binary relation$\bar{\chi}$ by
$$
((\alpha,a),\bar{z})\in\bar{\chi}\,\,\Leftrightarrow\,\,
(\bar{z},a)\in\bar{\varphi}_\alpha
$$
we analogously obtain the projection 
$\bar{\pi}=\bar{\vartheta}^{-1}\circ\bar{\chi}$, where
$\bar{\vartheta}$ is a bijection of the canonical 
decomposition of surjection $\bar{\chi}$.

Since by $(\ref{4})$ we have $\chi^{-1}\circ\chi=
\bar{\chi}^{-1}\circ\bar{\chi}$ then $\pi=\bar{\pi}$ and therefore
$\chi=\vartheta\circ \bar{\vartheta}^{-1}\circ\bar{\chi}$. 
Definitions of relation $\chi$ and relation $\bar{\chi}$ give $(\ref{5})$,
where
$$
\omega=\bar{\vartheta}\circ\vartheta^{-1}\colon \cup_{\beta\in
I}\dom\varphi_\beta\to\cup_{\beta\in I}\dom\bar{\varphi}_\beta
$$
is a bijection that is an injective co-injection if considered 
as a binary relation. $\Box$

\section*{Connections to solution of Sincov's equation on a group}

Given the composition operation of binary relations, 
the codomain $\cod\Phi=2^{M\times M}$ 
has a structure of a monoid with identity element $\id_M$. 
Replacing in $(\ref{1})$ all inequalities by equalities, 
the image $\img\Phi$ with regard to the same operation 
has a group structure because
for any $\Phi_{\alpha\beta}\in\img\Phi$ there exists an inverse 
element which is  coincidentally the inverse relation
$\Phi_{\alpha\beta}^{-1}=\Phi_{\beta\alpha}\in\img\Phi$.
If we fix in $(\ref{1})$ the index value $\gamma\in I$ and we denote 
$\varphi_\alpha=\Phi_{\alpha\gamma}$, we obtain 
$\Phi_{\alpha\beta}\circ\varphi_\beta=\varphi_\alpha$. Hence, applying 
group properties it follows that 
$$
\Phi_{\alpha\beta}=\varphi_\alpha\circ\varphi_\beta^{-1}
$$
and this is the general solution of Sincov's equation
$\Phi_{\alpha\beta}\circ\Phi_{\beta\gamma}=\Phi_{\alpha\gamma}$
in case $\circ$ is a group operation. This solution has been already 
published and described in detail, e.g., in Theorem 2 in
8.1.3 \cite{5} and has used only the fact that associativity and 
invertibility hold for any group and not only for Abelian groups 
as originally considered by Sincov \cite{6}.

Equality $(\ref{1})$ gives $\Phi_{\alpha\beta}\circ
\Phi_{\alpha\beta}^{-1}=\Phi_{\alpha\beta}^{-1}\circ
\Phi_{\alpha\beta}=\id_M$, which means that $\Phi_{\alpha\beta}$ and
then consequently also $\varphi_\alpha$ are bijections $M\to M$. 
As such they are also injective co-injections, hence the general solution of 
the modified Sincov's equation 8.1.3(9) \cite{5} is a particular
case of our general solution $(\ref{3})$ of inequality $(\ref{1})$.

\section*{Differentiable manifolds}

Let $p\geq0$ be an integer, let $M$ be a Banach space. Since we 
can consider the values $\Phi_{\alpha\beta}$ of solution $\Phi$ of 
functional inequalities $(\ref{1})$ as maps, we can impose 
differentiability requirements on them.  If we insist all the values 
$\Phi_{\alpha\beta}$ of general solution $(\ref{3})$ are
 $C^p$--isomorphisms of open sets, then such solution has a 
structure of a differentiable manifold of class $C^p$. In particular, 
the system of pairs $\{(\dom\varphi_\alpha, \varphi_\alpha)\}_{\alpha\in I}$ 
is therefore by  $(\ref{2})$ and $(\ref{3})$ a $C^p$--{\it atlas\/} 
on the set $X=\cup_{\alpha\in I}\dom\varphi_\alpha$ in the sense of
Lang's axioms  \cite{7}, Chapter II, \S1:

\medskip

{\bf AT1.} {\it Each $U_\alpha$ is a subset of 
$X$ and $\{U_\alpha\}_{\alpha\in I}$ cover $X$.

{\bf AT2.} Each $\varphi_\alpha$ is a bijection of 
$U_\alpha$ onto an open subset $\varphi_\alpha(U_\alpha)$
of the Banach space $M$ and for any $\alpha,\beta$, 
$\varphi_\alpha(U_\alpha\cap U_\beta)$ is open in $M$.

{\bf AT3.} The map
$$
\Phi_{\alpha\beta} \colon \varphi_\beta(U_\beta\cap U_\alpha) \ni
a\mapsto \varphi_\alpha(\varphi_\beta^{-1}(a))\in
\varphi_\alpha(U_\alpha\cap U_\beta)
$$
is a $C^p$--isomorphism for each pair of indices $\alpha,\beta$,\/}

\medskip

\noindent where $U_\alpha$ is another notation for 
$\dom\varphi_\alpha$. The other way around, any $C^p$--atlas 
satisfying these axioms generates by $(\ref{3})$ solution
of functional inequalities $(\ref{1})$.

The set $X=\cup_{\alpha\in I}\dom\varphi_\alpha$ equipped with 
a $C^p$--atlas $\{(\dom\varphi_\alpha, \varphi_\alpha)\}_{\alpha\in
I}$ is called a {\it differentiable manifold of class $C^p$\/} or 
shortly a  {\it $C^p$--manifold\/},
and the maps $\Phi_{\alpha\beta}$ appearing in the axiom {\bf AT3} 
are called  {\it transition maps\/} (see, e.g., Chapter 1 in
\cite{8}). The co-injection $\omega\colon X\to\bar{X}$
in Theorem 2 is therefore an {\it isomorphism of $C^p$--manifolds\/}.

An example of a $C^1$--manifold is the manifold of solutions of 
ordinary differentiable equation described in our motivation 
example above.

\bigskip

\subsection*{Remark on Hausdorff property.} From the above cited 
Lang's a\-xi\-o\-ms as they appear in the first edition of its book 
\cite{7}, the Hausdorff condition of a suitable differentiable manifold 
$X$ does not follow. But thanks to axiom {\bf AT3} we get
$$
\Phi_{\alpha\beta}\subseteq\rng\Phi_{\beta\beta}\times\rng\Phi_{\alpha\alpha}
$$
and it is easy to see, that $X$ is a {\it Hausdorff space\/} if and only if
for any pair $(\alpha,\beta)\in I\times I$ the set
$\Phi_{\alpha\beta}$ is closed in 
$\rng\Phi_{\beta\beta}\times\rng\Phi_{\alpha\alpha}$.
The Hausdorff property of the manifold can be achieved in a similar
way as its differentiability by means of demands placed only on 
transition maps.

\section*{Sincov's inequality for transition relations}

From the above it follows that if the solution values $\Phi_{\alpha\beta}$ $(\ref{2})$
of inequalities $(\ref{1})$ are $C^p$--isomorphisms of open sets, then we can 
interpret them as transition maps and the system of injective
co-injections $\{\varphi_\alpha\}_{\alpha\in I}$ from Theorem 1 as 
their corresponding atlas on the set $X=\cup_{\alpha\in I}\dom\varphi_\alpha$.

The above interpretation makes sense even though $\Phi_{\alpha\beta}$
are not $C^p$-isomorphisms of open sets. Hence, we are able to 
generalize the notion of system of transition maps to the term 
{\it system of transition relations\/}, that we define as a system of binary
relations $\{\Phi_{\alpha\beta}\}_{(\alpha,\beta)\in I\times I}$
solving inequalities $(\ref{1})$ representing therefore {\it
Sincov's inequality for transition relations\/}. By Theorem 1, for any 
system of transition relations there exists an {\it atlas\/}
$\{\varphi_\alpha\}_{\alpha\in I}$ which is by Theorem 2 unique up to 
some bijection $\omega$.

Consequently, Sincov's inequality for transition relations represents 
a generalization of the theory on differentiable manifolds to 
non-differen\-tiable and discontinuous case. This allows to pose 
many important questions for future work.  A direct application is 
expected for generalization in the analysis of our motivation example 
from differential to difference equations.

\end{document}